\documentclass[12pt]{amsart}

\usepackage[T1]{fontenc}
\usepackage{lmodern}
\usepackage[letterpaper,margin=1in]{geometry}
\usepackage{amssymb,mathtools,microtype}
\usepackage{enumitem,booktabs}
\usepackage[hidelinks]{hyperref}
\hypersetup{
  pdftitle={Depth of Powers of Cover Ideals of Uniform Path Hypergraphs},
  pdfsubject={Cover ideals, linear quotients, and depth of powers},
  pdfkeywords={cover ideal, path hypergraph, symbolic power, linear quotients, depth}
}

\newtheorem{theorem}{Theorem}[section]
\newtheorem{proposition}[theorem]{Proposition}
\newtheorem{lemma}[theorem]{Lemma}
\newtheorem{corollary}[theorem]{Corollary}
\theoremstyle{definition}
\newtheorem{definition}[theorem]{Definition}
\newtheorem{example}[theorem]{Example}
\theoremstyle{remark}

\numberwithin{equation}{section}

\newcommand{\avec}{\mathbf a}

\newcommand{\xvec}{\mathbf x}
\newcommand{\GG}{\mathcal G}

\DeclareMathOperator{\depth}{depth}
\DeclareMathOperator{\pd}{pd}
\DeclareMathOperator{\reg}{reg}

\setlist[enumerate]{label=\textup{(\roman*)},leftmargin=2em,itemsep=3pt}
\title[Depth of powers of cover ideals]{Depth of Powers of Cover Ideals of Path Ideals}

\author[N.T. Hang]{Nguyen Thu Hang}
\address{Thai Nguyen University of Sciences, Phan Dinh phung Ward, Thai Nguyen, Vietnam}
\email{hangnt@tnus.edu.vn}

\author{Truong Thi Hien}
\address{Hong Duc University, 565 Quang Trung Street, Hac Thanh Ward, Thanh Hoa, Vietnam}
\email{hientruong86@gmail.com}

\author{Thanh Vu}
\address{Institute of Mathematics, VAST, 18 Hoang Quoc Viet, Hanoi, Vietnam}
\email{vuqthanh@gmail.com}

\date{}

\begin{document}

\begin{abstract}
Let $ J_t(P_n) = (x_1, \ldots, x_t) \cap \cdots \cap (x_{n-t+1}, \ldots, x_n) $ be the cover ideal of the $t$-path ideal of a path graph inside the polynomial ring $ S = K[x_1, \ldots, x_n] $ over a field $K$. For every $ n \ge t \ge 2 $ and $s \ge 2$, we prove that
\[
\depth \left ( S/J_t(P_n)^s \right ) = q(s+1) + \left\lfloor \frac{r-1}{t} \right\rfloor,
\]
where $ n = q(ts + 1) + r $ with $ 1 \le r \le ts + 1 $. Our proof shows that all powers of $J_t(P_n)$ have linear quotients in increasing lexicographic order and explicitly describes the corresponding colon ideals.
\end{abstract}

\maketitle

\section{Introduction}\label{sec:introduction}
By the foundational result of Brodmann \cite{B}, the depth of powers eventually becomes constant. Nonetheless, the behavior of the depth of powers is poorly understood, even for edge ideals of graphs. Recently, computational efforts have increasingly focused on powers of cover ideals. Notably, the depth of cover ideals of graphs has been determined for paths \cite{DHNT}, cycles \cite{DHV}, and Ferrers graphs \cite{HV}. In this direction, the techniques introduced by Dung, Hang, and Vu \cite{DHV} for computing the depth of symbolic powers of cover ideals of hypergraphs-by leveraging the regularity of admissible subhypergraphs-have played a vital role. Outside the class of edge ideals, however, the depth of powers of cover ideals of general hypergraphs remains largely unknown. In this work, we compute the depth of powers of cover ideals of $t$-path ideals of path graphs. To the best of our knowledge, these are the first non-trivial results for general $t$-path graphs. Although the techniques of Dung, Hang, and Vu \cite{DHV} can be adapted to this setting-noting that the subhypergraphs of $t$-path hypergraphs of paths are simplicial forests, which allows one to apply the results of Erey and Faridi \cite{EF}-we take a different approach here. We show that all powers of the cover ideals of $t$-path ideals of paths admit linear quotients, a notion introduced by Herzog and Takayama \cite{HT}, in increasing lexicographic order, and we explicitly describe the corresponding colon ideals. We note that Bijender, Kumar, and Kumar \cite{BKK} have shown more generally that powers of cover ideals of simplicial trees have linear quotients; however, their results do not provide an explicit description of the quotients. Our explicit linear order allows us to describe these quotients precisely. We now introduce the necessary notation and describe our main results in detail.

Let $K$ be a field and let $S = K[x_1, \ldots, x_n]$ be the standard graded polynomial ring. Fix integers $n \ge t \ge 2$, and denote the index set by $[n] = \{1, \ldots, n\}$. For each $i \in [n-t+1]$, let
\[
e_i = \{i, i+1, \ldots, i+t-1\}.
\]
The $t$-path ideal of $P_n$ and its cover ideal are given by
$$I_t(P_n) = (x_1 x_2 \cdots x_t, \ldots, x_{n-t+1} \cdots x_n) \quad \text{and} \quad J_t(P_n) = \bigcap_{i=1}^{n-t+1} (x_i, x_{i+1}, \ldots, x_{i+t-1}).$$
For brevity, set $J=J_t(P_n)$.

For integers $u, v$, write $[u, v] = \{u, u+1, \ldots, v\}$; this set is empty if $u > v$. We use $\mathbb{N}$ for the set of nonnegative integers. For $\mathbf{a} \in \mathbb{N}^n$, set
\[
\mathbf{x}^{\mathbf{a}} = x_1^{a_1} \cdots x_n^{a_n}, 
\qquad |\mathbf{a}| = a_1 + \cdots + a_n, 
\qquad w_i(\mathbf{a}) = \sum_{j=i}^{i+t-1} a_j.
\]
The minimal monomial generating set of a monomial ideal $I$ is denoted by $\mathcal{G}(I)$. By \cite[Theorem 1.1]{HHT}, $J_t(P_n)^{(s)} = J_t(P_n)^s$ for every $s \ge 1$, where $J_t(P_n)^{(s)}$ is the $s$-th symbolic power, which is the ideal generated by the $s$-covers of $I_t(P_n)$. Formally, we have:

\begin{definition}
An \emph{integral $s$-cover} of $I_t(P_n)$ is a vector $\mathbf{a} \in \mathbb{N}^n$ satisfying $w_i(\mathbf{a}) \ge s$ for every $1 \le i \le n-t+1$. It is \emph{minimal} if there is no distinct integral $s$-cover $\mathbf{b}$ such that $\mathbf{b} \le \mathbf{a}$ componentwise. For a fixed integral $s$-cover $\mathbf{a}$, the interval $e_i$ is \emph{tight} if $w_i(\mathbf{a}) = s$.
\end{definition}

Consequently,
$$
\mathbf{x}^{\mathbf{a}} \in J_t(P_n)^s \quad \Longleftrightarrow \quad w_i(\mathbf{a}) \ge s \quad \text{for all } 1 \le i \le n-t+1.$$

We now describe the linear order and the quotients. 
Fix $s \ge 1$ and $\mathbf{x}^{\mathbf{a}} \in \mathcal{G}(J_t(P_n)^s)$.
For every $k$ with $a_k > 0$, define
\begin{align*}
\lambda_{\mathbf{a}}(k) &= \min\{i : k \in e_i, \ w_i(\mathbf{a}) = s\}, \\
\rho_{\mathbf{a}}(k) &= \lambda_{\mathbf{a}}(k) + t - 1, \\
R_{\mathbf{a}}(k) &= [k+1, \rho_{\mathbf{a}}(k)].
\end{align*}
The \emph{exchange set} associated with $\mathbf{a}$ is
$$
Q(\mathbf{a}) = \bigcup_{a_k > 0} R_{\mathbf{a}}(k).
$$
In particular, $1 \notin Q(\mathbf{a})$.

We order exponent vectors increasingly in lexicographic order: $\mathbf{b} \prec \mathbf{a}$ if, at the first index $k$ at which they differ, $b_k < a_k$. The same notation is used for the associated monomials. Recall that a monomial ideal has \emph{linear quotients} with respect to an ordering $u_1, \ldots, u_N$ of its minimal generators if every colon ideal $(u_1, \ldots, u_{r-1}) : u_r$ for $r \ge 2$ is generated by variables.

\begin{theorem}\label{thm:linear-quotients}
Let $u_1 \prec \cdots \prec u_N$ be the minimal monomial generators of $J_t(P_n)^s$ in increasing lexicographic order. If $r > 1$ and $u_r = \mathbf{x}^{\mathbf{a}}$, then
$$
(u_1, \ldots, u_{r-1}) : u_r = (x_j : j \in Q(\mathbf{a})).$$
In particular, $J_t(P_n)^s$ has linear quotients for every $s \ge 1$.
\end{theorem}

Hence, to determine the projective dimension of $J_t(P_n)^s$, it suffices to find $\max_{\mathbf{a}} |Q(\mathbf{a})|$ over all exponents $\mathbf{a}$ corresponding to the minimal generators of $J_t(P_n)^s$. This leads us to our main theorem.

\begin{theorem}\label{thm:main}
Let $n \ge t \ge 2$ and $s \ge 2$. Set $L = ts + 1$, and write
\[
n = qL + r, \qquad q \ge 0, \qquad 1 \le r \le L.
\]
Then
$$
\depth_S(S / J_t(P_n)^s) = q(s+1) + \left\lfloor \frac{r-1}{t} \right\rfloor.$$
\end{theorem}

In Section~\ref{sec:covers}, we describe the minimal generators of $J_t(P_n)^s$ and the linear ordering that yields its linear quotients. In Section~\ref{sec:bound}, we establish the upper bound for the cardinality of the exchange sets. Finally, in Section~\ref{sec:extremal}, we construct a minimal generator attaining this bound and deduce Theorem~\ref{thm:main}.

\section{Minimal generators and linear quotients}\label{sec:covers}
Throughout, $S = K[x_1, \ldots, x_n]$ is a standard graded polynomial ring over a field $K$. For a finitely generated graded $S$-module $M$, we write
\[
\beta_{i,j}(M) = \dim_K \operatorname{Tor}_i^S(M, K)_j
\]
for its graded Betti numbers.

\begin{definition}
Let $M$ be a finitely generated graded $S$-module. The projective dimension and the Castelnuovo--Mumford regularity of $M$ are defined by
\begin{align*}
\pd_S(M) &= \max \{ i \mid \beta_{i,j}(M) \neq 0 \text{ for some } j \}, \\
\reg_S(M) &= \max \{ j - i \mid \beta_{i,j}(M) \neq 0 \text{ for some } i, j \}.
\end{align*}
\end{definition}

First, we give a description of the minimal generators of $J_t(P_n)^s$.

\begin{lemma}\label{lem:minimality}
Let $\mathbf{a} \in \mathbb{N}^n$ be an integral $s$-cover. Then $\mathbf{x}^{\mathbf{a}} \in \mathcal{G}(J_t(P_n)^s)$ if and only if every index $k$ with $a_k > 0$ belongs to a tight interval. In particular, every minimal integral $s$-cover satisfies $0 \le a_k \le s$.
\end{lemma}

\begin{proof}
Suppose $\mathbf{x}^{\mathbf{a}} \in \mathcal{G}(J_t(P_n)^s)$ and $a_k > 0$. If $w_i(\mathbf{a}) \ge s+1$ for every $i$ with $k \in e_i$, then decreasing $a_k$ by $1$ would still yield an integral $s$-cover, contradicting the minimality of $\mathbf{x}^{\mathbf{a}}$. Conversely, suppose that for every $k$ with $a_k > 0$ there exists an index $i$ such that $k \in e_i$ and $w_i(\mathbf{a}) = s$. Then decreasing any positive coordinate $a_k$ would decrease $w_i(\mathbf{a})$ below $s$, so the resulting monomial would no longer lie in $J_t(P_n)^s$. Hence, $\mathbf{x}^{\mathbf{a}}$ is a minimal generator. 

Finally, for any $k \in [n]$, choose an index $i$ such that $k \in e_i$. If $a_k > s$, then $w_i(\mathbf{a}) \ge a_k \ge s+1$ for every such $i$, so $k$ cannot belong to any tight interval; thus $0 \le a_k \le s$.
\end{proof}

\begin{proof}[Proof of Theorem \ref{thm:linear-quotients}]
Let $L_{\mathbf{a}} = (x_j \mid j \in Q(\mathbf{a}))$. 
Let $j \in Q(\mathbf{a})$. Choose $k$ with $a_k > 0$ and $k < j \le \rho_{\mathbf{a}}(k)$, and set $\mathbf{b} = \mathbf{a} - \mathbf{e}_k + \mathbf{e}_j$, where $\mathbf{e}_j$ denotes the unit vector whose only non-zero entry is $1$ in the $j$-th position. Let $e_i$ be a tight interval containing $k$. By definition, $\lambda_{\mathbf{a}}(k) \le i  \le k$. Hence, $i \le k < j \le \rho_{\mathbf{a}}(k) = \lambda_{\mathbf{a}}(k) + t -1 \le i+t-1$. In other words, $j \in e_i$, which implies $w_i(\mathbf{b}) \ge s$. For every other interval containing $k$ that is not tight, we still have $w_i(\mathbf{b}) \ge s$. If an interval does not contain $k$, then clearly $w_i(\mathbf{b}) \ge s$ holds as well. Consequently, $\mathbf{b}$ is an integral $s$-cover.

Choose $v = \mathbf{x}^{\mathbf{c}} \in \mathcal{G}(J^s)$ dividing $\mathbf{x}^{\mathbf{b}}$. 
Since $c_i \le a_i$ for $i < k$ and $c_k < a_k$, one has $v \prec \mathbf{x}^{\mathbf{a}}$. The only coordinate in which $\mathbf{b}$ exceeds $\mathbf{a}$ is $j$, with an excess of one. Consequently,
\[
\frac{v}{\gcd(v, \mathbf{x}^{\mathbf{a}})} = x_j.
\]
Therefore, $L_{\mathbf{a}} \subseteq (u_1, \ldots, u_{r-1}) : u_r$, where $u_r = \mathbf{x}^{\mathbf{a}}$.

Conversely, let $v = \mathbf{x}^{\mathbf{b}} \prec \mathbf{x}^{\mathbf{a}}$ be an earlier minimal generator, and let $k$ be the first index at which $\mathbf{b}$ and $\mathbf{a}$ differ. Then $b_k < a_k$, so $a_k > 0$. On the tight interval $e_{\lambda_{\mathbf{a}}(k)}$, the sum of the coordinates of $\mathbf{a}$ is $s$, whereas the corresponding sum for $\mathbf{b}$ is at least $s$. Since $b_i = a_i$ for $i < k$, there must exist an index
\[
k < j \le \rho_{\mathbf{a}}(k) \quad \text{such that} \quad b_j > a_j.
\]
Thus, $j \in Q(\mathbf{a})$, and $x_j$ divides $v / \gcd(v, \mathbf{x}^{\mathbf{a}})$. Hence, $v / \gcd(v, \mathbf{x}^{\mathbf{a}}) \in L_{\mathbf{a}}$. The conclusion follows.
\end{proof}

\begin{corollary}\label{cor:betti}
For every $s \ge 1$ and $i, d \ge 0$,
\[
\beta_{i,i+d}(J^s) = \sum_{\substack{\xvec^\avec \in \GG(J^s) \\ |\avec| = d}} \binom{|Q(\avec)|}{i},
\]
where $\binom{q}{i} = 0$ for $i > q$. In particular,
\[
\pd_S(J^s) = \max_{\xvec^\avec \in \GG(J^s)} |Q(\avec)|.
\]
\end{corollary}
\begin{proof}
The first claim follows from Theorem \ref{thm:linear-quotients} and \cite[Corollary~2.7]{SV}. The formula for the projective dimension then follows directly from the definition.
\end{proof}
\section{An upper bound for the exchange sets}\label{sec:bound}

Fix $s \ge 1$ and $\mathbf{x}^{\mathbf{a}} \in \mathcal{G}(J^s)$. The goal of this section is to prove the desired upper bound for $|Q(\mathbf{a})|$. To do so, we analyze the complement of $Q(\mathbf{a})$ using the following notation: $\operatorname{supp}(\mathbf{a}) = \{i \mid a_i \neq 0\}$ and $p = \min \operatorname{supp}(\mathbf{a})$. 

Let $T = \{e \mid w_{e-t+1}(\mathbf{a}) = s\}$ be the set of right endpoints of the tight intervals of $\mathbf{a}$. First, we have the following lemma.

\begin{lemma}\label{lem:complements}
Both the first and last intervals, $[1, t]$ and $[n-t+1, n]$, are tight, and
\[
[1, n] \setminus Q(\mathbf{a}) = [1, p] \cup \{e + 1 \mid e \in T, e < n\}.
\]
\end{lemma}

\begin{proof}
Since $w_1(\mathbf{a}) \ge s$, we deduce that $p \le t$. Let $h = \lambda_{\mathbf{a}}(p)$ be the starting index of the first tight interval containing $p$; in other words, $w_h(\mathbf{a}) = s$. Furthermore, $h \le p$, which implies $a_1 = \cdots = a_{h-1} = 0$. Therefore, $w_1(\mathbf{a}) \le w_h(\mathbf{a}) = s$. Combined with $w_1(\mathbf{a}) \ge s$, this yields $w_1(\mathbf{a}) = s$, meaning the first interval $[1, t]$ is tight. By symmetry, the last interval $[n-t+1, n]$ is also tight. By definition, every element in $Q(\mathbf{a})$ is greater than $p$, hence $[1, p] \subseteq [1, n] \setminus Q(\mathbf{a})$.

Now, let $e \in T$ with $e < n$. For any index $k$ with $a_k > 0$ and $k \le e$, we have $\rho_{\mathbf{a}}(k) \le e$. Indeed, if $k \ge e-t+1$, then the definition of $T$ forces $\lambda_{\mathbf{a}}(k) \le e-t+1$, which implies $\rho_{\mathbf{a}}(k) \le e$ and hence $e + 1 \notin R_{\mathbf{a}}(k)$.

Conversely, let $j \notin Q(\mathbf{a})$ with $j > p$. Let $i$ be the largest index in $\operatorname{supp}(\mathbf{a})$ such that $i < j$, and let $e = \rho_{\mathbf{a}}(i)$. Since $j \notin Q(\mathbf{a})$, we must have $e < j$. By definition, $a_{i+1} = \cdots = a_{j-1} = 0$. Consequently, $w_{j-t}(\mathbf{a}) = a_{j-t} + \cdots + a_{j-1} \le w_{e-t+1}(\mathbf{a}) = s$. Since $\mathbf{a}$ is an integral $s$-cover, we have $w_{j-t}(\mathbf{a}) = s$, which implies $j - 1 \in T$. The conclusion follows.
\end{proof}

Assume that $[1, n] \setminus Q(\mathbf{a}) = \{b_1 < \ldots < b_m\}$, and set $b_{m+1} = n+1$ and set $g_j = b_{j+1} - b_j$ for $j = 1, \ldots,m$. Then we have the following gap inequality.

\begin{lemma}\label{lem:gap-inequality}
We have $1 \le g_j \le t$ for all $j = 1, \ldots, m$, and
\[
g_{i} + \cdots + g_{i+s} \le ts + 1
\]
for every $i = 1, \ldots, m-s$.
\end{lemma}

\begin{proof}
Let $\ell>k$ be two consecutive elements of $T$.
Suppose that $\ell-k>t$. The interval $e_{\ell-t}$
is not tight, whereas $e_{\ell-t+1}$ is tight. Hence
\[
w_{\ell-t}(\mathbf a)-w_{\ell-t+1}(\mathbf a)
=a_{\ell-t}-a_\ell\ge1.
\]
Thus $a_{\ell-t}>0$. By Lemma~\ref{lem:minimality},
some tight interval contains $\ell-t$. Its right
endpoint belongs to $[\ell-t,\ell-1]$, which is
contained in $(k,\ell)$, a contradiction. Therefore
$\ell-k\le t$. Since $t,n\in T$ and $p\le t$,
Lemma~\ref{lem:complements} gives
$1\le g_j\le t$ for every $j$.

Suppose that
$
g_i+\cdots+g_{i+s}\ge ts+2
$
for some $1\le i\le m-s$. Set
$
\delta_j=t-g_{i+j}\ (0\le j\le s).
$
Then
$
\delta_j\ge0,\
\sum_{j=0}^{s}\delta_j\le t-2.
$
In particular, $g_{i+j}\ge2$ for every $0\le j\le s$,
so $b_i\ge p$. Set
$
c_j=b_{i+j}-1\ (0\le j\le s+1).
$
By Lemma~\ref{lem:complements},
$c_1,\ldots,c_{s+1}$ are consecutive elements of $T$.
Moreover, $c_0\in T$ if $b_i>p$, and
$
c_{j+1}-c_j=t-\delta_j\ (0\le j\le s).
$

For $j = 0, \ldots, s$, define

\[h_0 = 1; \quad  
h_j = 1 + \sum _{\ell=0}^{j-1} \delta_\ell \quad \text{and} \quad P_j = \sum_{v = c_j + 1}^{c_j + h_j} a_v.
\]

We will show that $P_0\ge1$ and
$P_j\ge P_{j-1}+1$ for $1\le j\le s$.

First, $P_0=a_{b_i}$. If $b_i=p$, then $P_0\ge1$.
If $b_i>p$, then $c_0\in T$. Since $g_i\ge2$,
the interval ending at $c_0$ is tight, whereas the
interval ending at $c_0+1$ is not. Consequently,
\[
w_{c_0-t+2}(\mathbf a)-w_{c_0-t+1}(\mathbf a)
=a_{b_i}-a_{c_0-t+1}\ge1.
\]
Thus $P_0\ge1$ in either case.

Now fix $1\le j\le s$. Since
$
h_j+\delta_j
=1+\sum_{\ell=0}^{j}\delta_\ell
\le t-1,
$
we have
$
1\le h_j<t-\delta_j=g_{i+j}.
$
The interval ending at $c_j$ is tight, whereas
the interval ending at $c_j+h_j$ is not tight.
Therefore
\[
1\le w_{c_j-t+1+h_j}(\mathbf a) -w_{c_j-t+1}(\mathbf a) =
P_j-\sum_{v=c_j-t+1}^{c_j-t+h_j}a_v.
\]
On the other hand,
$
c_{j-1}+1
=c_j-t+1+\delta_{j-1}
\ge c_j-t+1, \text{and} \
c_{j-1}+h_{j-1}
=c_j-t+h_j.
$
Hence
$
P_{j-1}\le
\sum_{v=c_j-t+1}^{c_j-t+h_j}a_v,
$
and $P_j\ge P_{j-1}+1$. It follows that $P_s\ge s+1$.

Finally,
$
c_{s+1}-t+1=c_s+1-\delta_s\le c_s+1,
$
and
$
c_s+h_s<c_s+g_{i+s}=c_{s+1}.
$
Thus
$
[c_s+1,c_s+h_s]
\subseteq[c_{s+1}-t+1,c_{s+1}].
$
Since the interval ending at $c_{s+1}$ is tight,
$
P_s\le w_{c_{s+1}-t+1}(\mathbf a)=s,
$
a contradiction.
\end{proof}

\begin{proposition}\label{prop:exchange-bound}
Let $n \ge t \ge 2$ and $s \ge 1$. Write $n = q(ts+1) + r$ with $1 \le r \le ts+1$. Then every $\mathbf{x}^{\mathbf{a}} \in \mathcal{G}(J^s)$ satisfies
\[
|Q(\mathbf{a})| \le n - 1 - q(s+1) - \left\lfloor \frac{r-1}{t} \right\rfloor.
\]
\end{proposition}

\begin{proof}
Let $[1, n] \setminus Q(\mathbf{a}) = \{b_1, \ldots, b_m\}$. Denote $b_{m+1} = n+1$ and $g_j = b_{j+1} - b_j$ for $j = 1, \ldots, m$, so that $n = \sum_{j=1}^m g_j$. Write $m = A(s+1) + B$ with $0 \le B \le s$. By Lemma \ref{lem:gap-inequality}, we have
\[
g_{i(s+1) + 1} + \cdots + g_{(i+1)(s+1)} \le ts + 1
\]
for every $i = 0, \ldots, A-1$. Consequently,
\[
n = \sum_{j=1}^m g_j \le A(ts+1) + Bt.
\]
For each integer $u$, define
\[
F(u) = \left\lfloor \frac{u}{s+1} \right\rfloor (ts+1) + \left(u - (s+1)\left\lfloor \frac{u}{s+1} \right\rfloor \right) t,
\]
which is an increasing function of $u$. Set $k = \lfloor \frac{r-1}{t} \rfloor$. If $m \le q(s+1) + k$, then
\[
n \le F(m) \le F(q(s+1) + k) = q(ts+1) + kt < n,
\]
which is a contradiction. Hence, $m \ge q(s+1) + k + 1$. The conclusion follows since $|Q(\mathbf{a})| = n - m$.
\end{proof}

\section{Minimal generators attaining the bound}\label{sec:extremal}
Throughout this section, let $s \ge 2$ and $L = ts + 1$. We regard exponent vectors as words, where $0^h$ denotes a string of $h$ zeros, and $(\ell)$ denotes a word of length $1$ with value $\ell$. For finite words $U$ and $V$, let $U \mid V$ denote their concatenation. We now introduce a construction of basic words and the general construction, which we will justify later.

For $1 \le j \le s$, define the following basic word of length $t$:
\[
B_j = (j, 0^{t-2}, s-j).
\]
Define the word of length $L$ by
\[
U = B_1 \mid B_2 \mid \cdots \mid B_s \mid (s).
\]
For each $t \le M \le L$, write $M = kt + d$ with $0 \le d \le t-1$ and define
\[
V_M = \begin{cases} 
B_1 \mid \cdots \mid B_k & \text{if } d = 0, \\ 
B_1 \mid \cdots \mid B_k \mid (k) & \text{if } d = 1, \\ 
B_1 \mid \cdots \mid B_{k-1} \mid (k, 0^{d-2}, s-k, 0^{t-d}, s, 0^{d-1}) & \text{if } d \ge 2. 
\end{cases}
\]

We now construct words of length $L + c$ for $1 \le c < t$:
\[
V_{L+1} = (s, 0^{t-1}, s) \mid B_1 \mid \cdots \mid B_{s-2} \mid (s, 0^{t-1}, s),
\]
and for $2 \le c < t$,
\[
V_{L+c} = (s, 0^{t-1}, s) \mid B_1 \mid \cdots \mid B_{s-2} \mid (s-1, 0^{c-2}, 1, 0^{t-c}, s, 0^{c-1}).
\]

For the general construction, write $n = qL + r$ with $1 \le r \le L$. We define the word/exponent vector $\mathbf{a}$ by
\begin{equation}\label{eq:construction}
\mathbf{a} = \begin{cases} 
U^q \mid V_r & \text{if } t \le r \le L, \\ 
U^{q-1} \mid V_{L+r} & \text{if } 1 \le r < t. 
\end{cases}
\end{equation}
In the second case, $q \ge 1$ since $n \ge t$. Hence, $\mathbf{a}$ is a well-defined word of length $n$.

\begin{example}\label{ex:n12}
Let $t = 3$ and $n = 12$. For $s = 2$, one has $L = 7$, $q = 1$, and $r = 5$. Hence, $\mathbf{a} = U \mid V_5$, where
\[
U = (1, 0, 1, 2, 0, 0, 2) \quad \text{and} \quad V_5 = (1, 1, 0, 2, 0).
\]
\end{example}

To verify that our construction yields a minimal generator of $J_t(P_n)^s$ attaining the maximal bound for $|Q(\mathbf{a})|$, we first prove this property for the basic words. In what follows, we say that a word (or vector) is minimal if it is a minimal integral $s$-cover. Recall from the previous section that, for a word/exponent vector $V$, we define
\[
T(V) = \{i \mid w_{i-t+1}(V) = s\}
\]
as the set of right endpoints of the tight intervals.

\begin{lemma}\label{lem:word_1}
Let $M = kt + d$ with $0 \le d < t$ and $t \le M \le L$. Then $V_M$ is minimal and
\[
T(V_M) = \begin{cases}
\{t, 2t, \ldots, kt\} & \text{if } d = 0, \\
\{t, 2t, \ldots, kt\} \cup \{M\} & \text{if } d > 0.
\end{cases}
\]
\end{lemma}

\begin{proof}
Each block $B_j$ has weight $s$. Every window crossing the boundary between $B_j$ and $B_{j+1}$ has weight $(s - j) + (j+1) = s + 1$, so the case $d = 0$ follows. When $d = 1$, the additional window at the end has weight $(s - k) + k = s$, which establishes the case $d = 1$. The case $2 \le d < t$ likewise follows by direct verification on the final concatenated words. Moreover, the description of $T(V_M)$ shows that every positive coordinate belongs to a tight window; hence, $V_M$ is minimal.
\end{proof}

\begin{lemma}\label{lem:words_2}
Let $1 \le c < t$. The word $V_{L+c}$ is minimal, and
\[
T(V_{L+c}) = \{t\} \cup \{jt + 1 \mid 1 \le j \le s\} \cup \{L + c\}.
\]
\end{lemma}

\begin{proof}
The first two windows ending at $t$ and $t+1$ both have weight $s$. Every window crossing from these blocks into the subsequent block contains the final entry $s$ along with a positive entry from the next block, and hence has weight strictly greater than $s$. The intermediate blocks $B_1, \ldots, B_{s-2}$ have weight $s$, and any window crossing them has weight greater than $s$. The final blocks can be verified directly to show that every window of length $t$ has weight at least $s$, yielding the described set $T(V_{L+c})$. This description also shows that every positive coordinate belongs to a tight window; hence, $V_{L+c}$ is minimal.
\end{proof}

\begin{lemma}\label{lem:construction_minimality}
The vector $\mathbf{a}$ defined in Equation \eqref{eq:construction} satisfies the condition that $\mathbf{x}^{\mathbf{a}}$ is a minimal generator of $J_t(P_n)^s$, and
\[
|T(\mathbf{a})| = q(s+1) + \left\lfloor \frac{r-1}{t} \right\rfloor + 1.
\]
\end{lemma}

\begin{proof}
Every copy of $U$ ends with the entry $s$, and every subsequent subword begins with a positive entry. A window of length $t$ crossing such a boundary contains both entries, so its weight is strictly larger than $s$. Consequently, every window of length $t$ of $\mathbf{a}$ has weight at least $s$, and every positive coordinate belongs to a tight window (which holds for each constituent block). By Lemma \ref{lem:minimality}, $\mathbf{x}^{\mathbf{a}}$ is a minimal generator of $J_t(P_n)^s$. 

Furthermore, by Lemma \ref{lem:word_1} and Lemma \ref{lem:words_2}, we deduce that
\[
|T(\mathbf{a})| = \begin{cases}
q(s+1) + \left\lfloor \frac{r-1}{t} \right\rfloor + 1 & \text{if } r \ge t, \\
(q-1)(s+1) + (s+2) = q(s+1) + 1 & \text{if } r < t.
\end{cases}
\]
The conclusion follows.
\end{proof}

\begin{proof}[Proof of Theorem~\ref{thm:main}]
Proposition~\ref{prop:exchange-bound} gives
\[
|Q(\mathbf{a})| \le n - 1 - q(s+1) - \left\lfloor \frac{r-1}{t} \right\rfloor
\]
for every $\mathbf{x}^{\mathbf{a}} \in \mathcal{G}(J^s)$, and Lemma~\ref{lem:construction_minimality} produces a minimal generator for which equality holds, since $|Q(\mathbf{a})| = n - |T(\mathbf{a})|$. The conclusion then follows from Corollary~\ref{cor:betti} and the Auslander--Buchsbaum formula.
\end{proof}

For a proper homogeneous ideal $I\subset S$, denote by
$\operatorname{dstab}(I)$ the least positive integer $s_0$ such that
$\depth(S/I^s)$ is constant for all $s\ge s_0$.

\begin{corollary}\label{cor:depth-stability}
Let $n\ge t\ge2$. The sequence
$\depth(S/J_t(P_n)^s)$ is nonincreasing for $s\ge2$, and
\[
\lim_{s\to\infty}\depth(S/J_t(P_n)^s)
=
\left\lfloor\frac{n-1}{t}\right\rfloor.
\]
Moreover,
\[
\operatorname{dstab}(J_t(P_n))
\le
\begin{cases}
\displaystyle
\max\left\{2,\left\lceil\frac{n-1}{2t}\right\rceil\right\},
& \text{if }t\mid(n-1),\\[6pt]
\displaystyle
\max\left\{2,\left\lceil\frac{n-1}{t}\right\rceil\right\},
& \text{if }t\nmid(n-1).
\end{cases}
\]
If the bound on the right is greater than $2$, then equality holds.
\end{corollary}

\begin{proof}
Fix $s\ge2$ and set
$q=\lfloor(n-1)/(ts+1)\rfloor$.
Theorem~\ref{thm:main} yields
\[
\depth(S/J_t(P_n)^s)
=
q+\left\lfloor\frac{n-1-q}{t}\right\rfloor
=
\left\lfloor\frac{n-1+(t-1)q}{t}\right\rfloor.
\]
Since $q$ is nonincreasing as $s$ increases and is eventually zero,
the first two assertions follow.

Write $n-1=at+b$, where $0\le b<t$. Then
\[
\depth(S/J_t(P_n)^s)
=
a+\left\lfloor\frac{b+(t-1)q}{t}\right\rfloor,
\]
and for $s\ge 2$,
\begin{equation}\label{eq:attain}
\depth(S/J_t(P_n)^s)=a
\quad\Longleftrightarrow\quad
b+(t-1)q<t .
\end{equation}
Since $t\ge2$, this is equivalent to
$q\le1$ if $b=0$, and to $q=0$ if $b>0$.

If $b=0$, then
\[
q\le1
\quad\Longleftrightarrow\quad
at<2(ts+1)
\quad\Longleftrightarrow\quad
a\le2s,
\]
where the last equivalence follows from $t\ge2$ and the integrality
of $a-2s$. Hence the limiting value is attained when
$s\ge\lceil(n-1)/(2t)\rceil$.
If $b>0$, then
\[
q=0
\quad\Longleftrightarrow\quad
n-1<ts+1
\quad\Longleftrightarrow\quad
s\ge\left\lceil\frac{n-1}{t}\right\rceil.
\]
Set
\[
s_*=
\begin{cases}
\left\lceil\dfrac{n-1}{2t}\right\rceil & \text{if } t\mid (n-1),\\[6pt]
\left\lceil\dfrac{n-1}{t}\right\rceil & \text{if } t\nmid (n-1).
\end{cases}
\]

For $s\ge 2$ we have $\depth(S/J_t(P_n)^s)=a$ if and only if $s\ge s_*$.
Thus $\depth(S/J_t(P_n)^s)=a$ for all $s\ge \max\{2,s_*\}$, which gives
$\operatorname{dstab}(J_t(P_n))\le\max\{2,s_*\}$.

If $s_*>2$, then $s_*-1\ge2$, and the preceding
equivalence gives
\[
\depth(S/J_t(P_n)^{s_*-1})>a,
\qquad
\depth(S/J_t(P_n)^{s_*})=a.
\]
Therefore
$\operatorname{dstab}(J_t(P_n))\ge s_*$, and equality
follows from the upper bound.
\end{proof}

\begin{example}\label{ex:s1}
The hypothesis $s \ge 2$ in Theorem~\ref{thm:main} cannot be omitted. For $t = 3$ and $n = 6$, the minimal generators of $J_3(P_6)$, listed in increasing lexicographic order, are
\[
x_3x_6, \quad x_3x_5, \quad x_3x_4, \quad x_2x_5, \quad x_2x_4, \quad x_1x_4.
\]
Their colon ideals are, respectively,
\[
(0), \quad (x_6), \quad (x_5, x_6), \quad (x_3), \quad (x_3, x_5), \quad (x_2, x_3).
\]
Consequently, $\operatorname{pd}_S(J_3(P_6)) = 2$ and $\operatorname{depth}_S(S/J_3(P_6)) = 3$, whereas the formula in Theorem~\ref{thm:main} (when formally applied for $s=1$) would instead give $2$ (since $6 = 1 \cdot 4 + 2$).
\end{example}

\noindent\textbf{Acknowledgment.} This work is partially supported by the Vietnam National Foundation for Science and Technology Development under the grant number NAFOSTED/101/10104/2026/1279.

\noindent\textbf{Declaration of Competing Interest.} The authors declare that they have no known competing financial interests or personal relationships that could have appeared to influence the work reported in this paper.


\begin{thebibliography}{9}

\bibitem[BKK]{BKK} Bijender, A. Kumar, and R. Kumar, {\it Powers of vertex cover ideals of simplicial trees}, Math. Nachr. {\bf 297}, (2024), 1116--1135.


\bibitem [B]{B} M. Brodmann, {\em The asymptotic nature of the analytic spread}, Math. Proc. Cambridge Philos. Soc. {\bf 86} (1979), 35--39.

\bibitem[DHNT]{DHNT}
Tran Duc Dung, Nguyen Thu Hang, Pham Hong Nam, and Nguyen Thi Thanh Tam,
\emph{The depth function of powers of cover ideals of path graphs},
preprint (2026), arXiv:2605.03347.
\url{https://arxiv.org/abs/2605.03347}.


\bibitem [DHV]{DHV} T. D. Dung, N. T. Hang, T. Vu, Admissible subgraphs and the depth of symbolic powers of cover ideals of graphs, arXiv:2605.03369v1.


\bibitem[EF]{EF}
N. Erey and S. Faridi,
\emph{Betti numbers of monomial ideals via facet covers}, Journal of Pure and Applied Algebra {\bf 220} (2016), 1990--2000.

\bibitem[HHT]{HHT}
J. Herzog, T. Hibi, and N. V. Trung,
\emph{Vertex cover algebras of unimodular hypergraphs},
Proc. Amer. Math. Soc. \textbf{137} (2009), no.~2, 409--414.
\href{https://doi.org/10.1090/S0002-9939-08-09308-8}
{doi:10.1090/S0002-9939-08-09308-8}.

\bibitem[HT]{HT} J. Herzog  and Y. Takayama, {\it Resolutions by mapping cones}, Homology Homotopy Appl. {\bf 4} (2002), 277--294.


\bibitem[HV]{HV} D. T. Hoang and T. Vu, {\em Projective dimension of powers of cover ideal of Ferrers graphs}, arXiv:2606.00772.

\bibitem[SV]{SV}
L. Sharifan and M. Varbaro,
\emph{Graded Betti numbers of ideals with linear quotients},
Le Matematiche \textbf{63} (2008), no.~2, 257--265.

\end{thebibliography}
\end{document}